\documentclass[11pt]{article}

\usepackage[colorlinks,citecolor=blue,urlcolor=blue]{hyperref}

\usepackage{amssymb}
\usepackage{amsmath}
\usepackage{amsthm}
\usepackage[authoryear]{natbib}

\newcommand{\PP}{{\rm P}}
\newcommand{\EE}{{\rm E}}

\newcommand{\lle}{\,\,{\lesssim}\,\,}

\newcommand{\eps}{\varepsilon}

\newcommand{\FF}{\mathcal{F}}
\newcommand{\GG}{\mathcal{G}}

\newcommand{\RR}{\mathbb{R}}
\newcommand{\NN}{\mathbb{N}}

\newcommand{\HHH}{\mathbb{H}}

\newcommand{\qv}[1]{\left< #1 \right>}

\newcommand{\BBB}{\mathbb{B}}

\renewcommand{\phi}{\varphi}

\newcommand{\given}{\,|\,}
\renewcommand{\l}{\lambda}

\newcommand{\e}{\varepsilon}

\newcommand{\norm}[1]{||#1||}

\newtheorem{theorem}{Theorem}[section]

\theoremstyle{definition}

\def\e{\varepsilon}

\def\b{\beta}

\def\l{\lambda}

\let\ge\geqslant

\let\le\leqslant
\def\bardelta{\overline{\delta}}

\begin{document}

\title{\bf \mbox{Optimality of  Poisson processes intensity} \mbox{learning with Gaussian processes}}

\author{
Alisa Kirichenko\thanks{Korteweg-de Vries Institute for Mathematics,
University of Amsterdam,
P.O. Box 94248, 1090 GE Amsterdam, The Netherlands. Email: a.kirichenko@uva.nl.}
 \and 
 Harry van Zanten\thanks{Korteweg-de Vries Institute for Mathematics,
University of Amsterdam,
P.O. Box 94248, 1090 GE Amsterdam, The Netherlands. Email: hvzanten@uva.nl.}
}

\date{March 2, 2015}

\maketitle

\begin{abstract}
In this paper we provide theoretical support for the so-called 
``Sigmoidal Gaussian Cox Process'' approach to learning the intensity 
of an inhomogeneous Poisson process on a $d$-dimensional domain. 
This method was proposed by Adams, Murray and MacKay (ICML, 2009), 
who developed 
a tractable computational approach and showed in simulation 
and real data experiments that it can work quite satisfactorily. 
The results presented in the present paper provide theoretical 
underpinning of the method. In particular, we show how to tune the priors 
on the hyper parameters of the model in order for the procedure to automatically 
adapt to the degree of smoothness of the unknown intensity and to achieve 
optimal convergence rates.

\bigskip

\noindent{\bf Running title:} Optimal Poisson intensity learning with GP's

\bigskip


\end{abstract}

\numberwithin{equation}{section}

\section{Introduction}

Inhomogeneous Poisson processes are widely used models for count and point data in a variety of applied areas.
A typical task in applications is to learn the underlying intensity of a Poisson process
from a realised point pattern. In this paper we consider nonparametric Bayesian approaches to this problem. 
These do not assume a specific parametric form of the intensity function and produce 
posterior distributions which do not only give an estimate of the intensity, e.g.\ through the 
posterior mean or mode, but also give a measure of the remaining uncertainty through the spread of the posterior.

Several papers have explored nonparametric Bayesian approaches in this setting. 
An early reference is \cite{moller}, who study log-Gaussian priors. 
\cite{Shota} recently considered Gaussian processes combined with different, non-smooth link functions.
Kernel mixtures priors are considered in \cite{kottas}. Spline-based priors are used in \cite{dimatteo} and \cite{serra}.

The present study is motivated by a method that is not covered by earlier theoretical papers, 
namely the method of \cite{comp_science}. These authors
presented the first approach that is also computationally fully nonparametric 
in the sense that it does not involve potentially inaccurate finite-dimensional approximations. 
The method involves a prior on the intensity that is a random multiple of a transformed Gaussian process (GP). 
Both the hyper parameters of the GP and the multiplicative constant are endowed with priors as well, 
resulting in a hierarchical Bayes procedure (details in Section \ref{sec: model}). 
Simulation experiments and real data examples in \cite{comp_science} show that 
the method can give very satisfactory results.

The aim of this paper is to advance the theoretical understanding of 
the method of \cite{comp_science}, 
which they termed ``Sigmoidal Gaussian Cox Process'' (SGCP). 
It is by now well known both from theory and practice that nonparametric Bayesian methods need to be tuned 
very carefully to produce good results. 
An unfortunate choice of the prior or incorrectly tuned hyper parameters can easily result in procedures 
that give misleading results or that make sub-optimal use of the information in the training data.
See for instance the by now classical reference \cite{DF}, or the more recent paper 
\cite{jlmr} and the references therein.

A challenge in this problem (and in nonparametric function learning in general) is to devise a procedure 
that avoids overfitting and underfitting. 
The difficulty is that the appropriate degree of ``smoothing'' depends on the (unknown) regularity 
of the intensity function that produces the data. Indeed, intuitively it is clear that if the function is 
very smooth then to learn the intensity at a certain location 
we can borrow more information from neighboring points than if it is very rough. 
Ideally we want to have a procedure that automatically uses the appropriate degree of smoothing, i.e.\
that {\em adapts} to regularity. 

To address this issue theoretically it is common to take an asymptotic point of view. 
Specifically, we assume that we have $n$ independent sets of training data, produced by Poisson processes on 
the $d$-dimensional domain $S = [0,1]^d$ (say), with the same intensity function $\lambda_0: S \to [0, \infty)$. We aim 
to construct the learning procedure such that we achieve an optimal learning rate, 
irrespective of the regularity level of the intensity.
In the problem at hand it is known that if $\lambda_0$ has regularity $\beta > 0$, then
the best rate that any procedure 
can achieve is of the order $n^{-\beta/(d+2\beta)}$. 
This can be made precise in the minimax framework, for instance. For a fixed estimation or learning 
procedure, one can determine the largest expected loss that is incurred when the true function 
generating the data 
is varied over a ball of functions with fixed regularity $\beta$, say. This will depend on $n$
and quantifies the worst-case rate of convergence for that fixed estimator for $\beta$-regular truths. 
The minimax rate
is obtained by minimising this over all possible estimators. So it is the best convergence 
rate that any procedure can achieve, uniformly over a ball of functions with fixed regularity $\beta$.
See, e.g., \cite{tsybakov} for a general introduction to the minimax approach and  \cite{kutoyants} or \cite{PTRF}
for minimax results in the context of the Poisson process model that we consider in this paper.

Note that the smoothness degree is unknown to us, so we can not use it in the construction of the procedure, 
but still we want that the posterior contracts around $\lambda_0$ at the 
 rate $n^{-\beta/(d+2\beta)}$, as $n \to \infty$, if $\lambda_0$ is $\beta$-smooth. 
In this paper we prove that with appropriate priors on the hyper parameters, the SGCP approach of 
\cite{comp_science} attains this optimal rate (up to a logarithmic factor). It does 
so for every regularity level $\beta > 0$, so it is fully {\em rate-adaptive}. 

Technically the paper uses the mathematical framework for studying contraction rates for Gaussian and 
conditionally Gaussian priors as developed in \cite{vdVvZGaussian} and \cite{vdVvZGamma}. 
 We also use an extended version of a general result for Bayesian inference for $1$-dimensional Poisson processes
from \cite{serra}. On a general level the line of reasoning is similar to that of \cite{vdVvZGamma}.
However, due to the presence of a link function and a random multiplicative constant in the SGCP model (see Section \ref{sec: model1}
ahead) the results of the latter paper do not apply in the present setting and additional mathematical arguments 
are required to prove the desired results.

The paper is organised as follows. In Section \ref{sec: model1} we describe the 
Poisson process observation model and the SGCP prior model, which together 
determine a full hierarchical Bayesian model. 
The main result about the performance of the SGCP approach is presented
and discussed in Section \ref{sec: main}. Mathematical proofs are
given in Section \ref{sec: proof}.
In Section \ref{sec: con} we make some concluding remarks.

\section{The SGCP model}
\label{sec: model1}

\subsection{Observation model}

We assume we observe $n$ independent copies of an inhomogeneous Poisson
process on the $d$-dimensional unit cube $S = [0,1]^d$ (adaptation to other domains is straightforward).
We denote these observed data by $N^1, \ldots, N^n$. 
Formally every $N^i$ is a counting measure on subsets of $S$.
The object of interest 
is the underlying {\em intensity function}. 
This is a (integrable) function $\lambda: [0,1]^d \to [0,\infty)$ with the property that given $\lambda$, 
every $N^j$ is a random counting measure on $[0,1]^d$ such that $N^j(A)$ and $N^j(B)$ are 
independent if the sets $A, B \subset [0,1]^d$ are disjoint and
 the number of points $N^j(B)$ falling in the set $B$ has a Poisson distribution with mean $\int_B\lambda(s)\,ds$. 
If we want to stress that the probabilities and expectations involving the observations $N^j$ depend
on $\lambda$, we use the notations $\PP_\lambda$ and $\EE_\lambda$, respectively. 
We note that instead of considering observations from $n$ independent Poisson processes with intensity $\l$, one could equivalently 
consider observations from a single Poisson process with intensity $n\lambda$.

\subsection{Prior model}

The SGCP model introduced in \cite{comp_science} postulates a-priori that the intensity 
function $\lambda$ is of the form 
\begin{equation}\label{eq: lambda}
\lambda(s) = \lambda^* \sigma(g(s)), \qquad s \in S, 
\end{equation}
where $\lambda^* > 0$ is an upper bound on $\lambda$, $g$ is a GP indexed by $S$ and 
$\sigma$ is the sigmoid, or {logistic} function on the real line, defined by 
$\sigma(x) = ({1+e^{-x}})^{-1}$.
In the computational section of \cite{comp_science} $g$ is modeled as a GP with squared 
exponential covariance kernel and zero mean, with a prior on the length scale parameter. 
The hyper parameter $\lambda^*$ is endowed with an independent gamma prior.

In the mathematical results presented in this paper we allow a bit more flexibility in the choice of the covariance kernel of the GP, 
the link function $\sigma$ and the priors on the hyper parameters. We assume that $g$ is a zero-mean, homogenous GP 
with covariance kernel given in spectral form by 
\begin{equation}\label{eq: g}
\EE g(s)g(t) = \int{e^{-i\left<\xi,\ell(t-s)\right>}\mu(\xi)\,d \xi}, \qquad s,t \in S, 
\end{equation}
where $\ell > 0$ is an (inverse) length scale parameter and $\mu$ is a spectral density on $\RR^d$ 
such that the map $a \mapsto \mu(a\xi)$ on $(0, \infty)$ is decreasing for every $\xi \in \RR^d$ and that satisfies 
\begin{equation}
\label{about_mu}
\int{e^{\delta ||\xi||}\mu (d \xi)} < \infty
\end{equation}
for some $\delta > 0$ (the Euclidean inner product and norm are denoted by $\qv{\cdot, \cdot}$ and $\|\cdot\|$, respectively).
Note that, in particular, the centered Gaussian spectral density satisfies this condition and corresponds to 
the squared exponential kernel 
\[
\EE g(s)g(t) = e^{{-\ell^2\|t-s\|^2}}.
\] 
We endow the length scale parameter $\ell$ with a prior with density $p_\ell$ on $[0, \infty)$, for which we assume the bounds,
for positive constants $C_1, D_1, C_2, D_2$, nonnegative constants $p$, $q$, and every sufficiently large $x > 0$,
\begin{equation}\label{eq: fl}
C_1 x^p \exp(-D_1 x^d \log^{q} {x}) \le p_\ell(x) \le C_2 x^p \exp(-D_2 x^d \log^{q} {x}).
\end{equation}
This condition is, for instance, satisfied if $\ell^d$ has a gamma distribution, which is a common choice in practice.
Note however that the technical condition \eqref{eq: fl} is only a condition on 
the tail of the prior on $\ell$. 
On the upper bound $\l^*$ we put a prior satisfying an exponential tail bound. Specifically, we use 
a positive, continuous prior density $p_{\l^*}$ on $[0, \infty)$ such that for some $c_0, C_0, \kappa > 0$, 
\begin{equation}\label{eq: flambda}
\int_{\lambda_0}^\infty p_{\l^*}(x) \,dx \le C_0 e^{-c_0\lambda_0^\kappa}
\end{equation}
for all $\lambda_0 > 0$. Note that this condition is fulfilled if we place a gamma prior on $\l^*$. 
Finally, we use a strictly increasing, infinitely smooth link function $\sigma: \RR \to (0,1)$ in \eqref{eq: lambda} that satisfies 
\begin{equation}\label{eq: sigma}
|\sqrt{\sigma(x)} - \sqrt{\sigma(y)}| \le c|x-y|
\end{equation}
for all $x, y \in \RR$. This condition is in particular fulfilled for the sigmoid function employed by \cite{comp_science}. 
It holds for other link functions as well, for instance for the cdf of the standard normal distribution.

\subsection{Full hierarchical model}
\label{sec: model}

With the assumptions made in the preceding section in place, 
the full hierarchical specification of the prior and observation model can then be summarised as follows:
\begin{align*}
\ell ~ & \sim p_\ell \quad \text{(satisfying \eqref{eq: fl})}\\
\lambda^* & \sim p_{\lambda^*} \quad \text{(satisfying \eqref{eq: flambda})}\\
g \given \ell, \lambda^* & \sim \text{GP with kernel given by \eqref{eq: g}--\eqref{about_mu}}\\
\lambda \given g, \ell, \lambda^* & \sim \text{defined by \eqref{eq: lambda}, 
with smooth $\sigma$ satisfying \eqref{eq: sigma}}\\
N^1, \ldots, N^n \given \lambda, g, \ell, \lambda^* & \sim 
\text{independent Poisson processes with intensity $\lambda$}.
\end{align*}
Note that under the prior, several quantities are, by construction, independent.
Specifically, $\ell$ and $\lambda_*$ are independent, and $g$ and $\lambda^*$ are 
independent.

The main results of the paper concern the posterior distribution of the intensity function 
$\lambda$, i.e.\ the conditional $\lambda \given N^1, \ldots, N^n$. 
Throughout we will denote the prior on $\lambda$ by $\Pi$ and the posterior 
by $\Pi(\cdot \given N^1, \ldots, N^n)$. 
In this setting Bayes' formula asserts that
\begin{equation}\label{eq: post}
\Pi(\l \in B \given N^1, \ldots, N^n) = \frac{\int_B p(N^1, \ldots, N^n\given \lambda) \,\Pi(d\lambda)}
{\int p(N^1, \ldots, N^n\given \lambda) \,\Pi(d\lambda)}, 
\end{equation}
where the likelihood is given by 
\[
p(N^1, \ldots, N^n\given \lambda) = \prod_{i=1}^n e^{\int_S \l(x) N^i(dx) - \int_S (\l(x)-1)\,dx} 
\]
(see, e.g., \cite{kutoyants}).

\section{Main result}
\label{sec: main}

Consider the prior and observations model described in the preceding section 
and let $\Pi(\cdot \given N^1, \ldots, N^n)$ be the corresponding posterior distribution of 
the intensity function $\l$. 

The following theorem describes how quickly the posterior distribution contracts 
around the true intensity $\lambda_0$ that generates the data. 
The rate of contraction depends on the smoothness level of $\l_0$. This is 
quantified by assuming that $\l_0$ belongs to the H\"older space $C^\beta[0,1]^d$ 
for $\beta > 0$. By definition a function on $[0,1]^d$ belongs to this space if it has 
 partial derivatives up to the order $\lfloor\beta\rfloor$ and if the $\lfloor\beta\rfloor$th
 order partial derivatives are all H\"older continuous of the order $\beta - \lfloor\beta\rfloor$.
 Here $\lfloor\beta\rfloor$ denotes the greatest integer strictly smaller than $\beta$. 
The rate of contraction is measured in the $L^2$-distance between the square root
of intensities. This is the natural statistical metric in this problem, as it can be 
shown that in this setting the Hellinger distance between the models with intensity functions $\l_1$
and $\l_2$ is equivalent to $\min\{\|\sqrt\l_1 - \sqrt\l_2\|_2, 1\}$ (see \cite{serra}).
Here $\|f\|_2$ denotes the $L^2$-norm of a function on $S = [0,1]^d$, i.e.\ 
$\|f\|^2_2 = \int_S f^2(s)\,ds$.

\begin{theorem}
\label{main_theorem}
Suppose that $\l_0\in C^{\b}[0,1]^d$ for some $\b>0$ and that $\l_0$ is strictly positive. Then for 
all sufficiently large $M>0$,
\begin{equation}\label{eq: main}
\EE_{\l_0}\Pi(\lambda: \|\sqrt{\lambda}-\sqrt{\lambda_0}\|_2 \ge M n^{-\beta/(d+2\beta)}\log^\rho n | N^1, \ldots, N^n) 
\to 0
\end{equation}
as $n \to \infty$, for some $\rho > 0$. 
\end{theorem}

The theorem asserts that if the intensity $\l_0$ that generates the data is $\beta$-smooth, then, asymptotically, 
all the posterior mass is concentrated in (Hellinger) balls around $\l_0$ with a radius 
that is up to a logarithmic factor of the optimal order $n^{-\beta/(d+2\beta)}$. 
Since the procedure does not use the knowledge of the smoothness level $\beta$, this indeed 
shows that the method is rate-adaptive, i.e.\ the rate of convergence 
adapts automatically to the degree of smoothness of the true intensity. 
Let us mention once again that the conditions of the theorem are in particular fulfilled if 
in \eqref{eq: lambda}, $\lambda^*$ is taken gamma, $\sigma$ is the sigmoid (logistic) function, 
and $g$ is a squared exponential GP with length scale $\ell$, with $\ell^d$ a gamma variable.

\section{Proof of Theorem \ref{main_theorem}}
\label{sec: proof}

To prove the theorem we employ an extended version of a result from \cite{serra} 
that gives sufficient conditions for having \eqref{eq: main} in the case $d=1$, cf.\ their Theorem 1. 
Adaptation to the case of a general $d \in \NN$ is straightforward. 
To state the result we need some (standard) notation and terminology. 
For a set of positive functions $\FF$ we write $\FF^c$ for its complement and $\sqrt{\FF}= \{\sqrt{f}, f\in\FF\}.$ For $\eps>0$ and a norm $\|\cdot\|$ on $\FF$, let $N(\eps, \FF, \norm{\cdot})$ be the minimal number of balls of radius $\eps$ with respect to norm $\|\cdot\|$ needed to cover $\FF$. 
The uniform norm $\|f\|_\infty$ of a function $f$ on $S$ is defined, as usual, as $\|f\|_\infty = \sup_{s \in S}|f(s)|$. 
The space of continuous function on $S$ is denoted by $C(S)$. 

Let $\Pi$ now be a general prior on the intensity function $\lambda$ and let $\Pi(\cdot \given N^1, \ldots, N^n)$ 
be the corresponding posterior \eqref{eq: post}.

\begin{theorem}
\label{serra_theorem}
Assume that $\lambda_0$ is bounded away from $0$. Suppose that for positive sequences $\overline{\delta}_n, {\delta}_n \to 0$ such that $n(\overline{\delta}_n \wedge {\delta}_n)^2 \to \infty$
\footnote{As usual, $a \wedge b = \min\{a,b\}$ and $a \vee b = \max\{a,b\}$.}
 as $n \to \infty$ and constants $c_1, c_2>0$, it holds that for all $L>1$, there exist subsets $\FF_n \subset C(S)$ and a constant $c_3$ such that
\begin{align}
&{1-\Pi(\FF_n) \le e^{-L n\delta^2_n}, \label{remaining_mass_0}} \\
&\Pi(\l: \norm{\lambda-\l_0}_\infty \le \delta_n) \ge c_1e^{-n c_2 \delta_n^2}, \label{prior_mass_0} \\
&\log N(\overline{\delta}_n, \sqrt{\FF_n}, \|\cdot\|_2) \le c_3 n\overline{\delta}^2_n \label{entropy_0}.
\end{align}
Then for $\eps_n=\overline{\delta}_n \vee \delta_n$ and all sufficiently large $M>0$,
\begin{equation}\label{thmcon}
\EE_{\l_0}\Pi(\lambda: \|\sqrt{\lambda}-\sqrt{\lambda_0}\|_2 \ge M \eps_n | N^1, \ldots N^n) 
\to 0
\end{equation}
as $n \to \infty.$
\end{theorem}

We note that this theorem has  a form that is commonly encountered in the literature 
on contraction rates for nonparametric Bayes procedures.  
The so-called ``prior mass condition'' \eqref{prior_mass_0} requires that the prior puts sufficient mass near
the true intensity function $\lambda_0$ generating the data. 
The ``remaining mass condition'' \eqref{remaining_mass_0} and the ``entropy condion'' \eqref{entropy_0}
together require that ``most'' of the prior mass should be concentrated on so-called ``sieves'' 
$\FF_n$ that are not too large in terms of their metric entropy. 
The sieves grow as $n \to \infty$ and in the limit they capture all the posterior mass.

In the subsequent subsections we will show that the prior 
defined in Section \ref{sec: model} fulfils the conditions of this theorem, for 
$\delta_n=n^{-{\beta}/({2\beta+d})}(\log{n})^{k_1}$ and $\bardelta_n= L_1 
n^{-{\beta}/({2\beta+d})}(\log n)^{({d+1})/{2}+2{k_1}}$, with $L_1 > 0$ and $k_1=({(1+d)\vee q})/({2+d/\beta})$.
The proofs build on earlier work, especially from \cite{vdVvZGamma}, in which results like 
\eqref{remaining_mass_0}--\eqref{entropy_0} have been derived for the GP's like $g$.  
Here we extend and adapt these results to deal with the additional link function $\sigma$
and the prior on the maximum intensity $\lambda^*$.

\subsection{Prior mass condition}

In this section we show that with $\lambda^*$, $\sigma$ and $g$ as specified in 
Section \ref{sec: model} and $\lambda_0 \in C^\beta(S)$, we have 
\begin{equation}\label{eq: tp1}
\PP(\|\lambda^*\sigma(g) - \lambda_0\|_\infty \le \delta_n) \ge c_1e^{-n c_2 \delta_n^2}
\end{equation}
for constants $c_1, c_2 > 0$ and $\delta_n$ as defined above.

The link function $\sigma$ is strictly increasing and smooth, hence it has 
a smooth inverse $\sigma^{-1}: (0,1) \to \RR$. Define the function $w_0$ on $S$ 
by 
\[
w_0(s) = \sigma^{-1}\Big(\frac{\lambda_0(s)}{2\|\lambda_0\|_\infty}\Big), \qquad s \in S,
\]
so that $\lambda_0 = 2\|\lambda_0\|_\infty\sigma(w_0)$. Since the function $\lambda_0$ 
is positive and continuous on the compact set $S$, it 
is 
bounded away from $0$ on $S$, say $\lambda_0 \ge a > 0$. 
It follows that ${\lambda_0(s)}/{2\|\lambda_0\|_\infty}$ varies in the compact 
interval $[{a}/{2\norm{\l_0}_{\infty}}, {1}/{2}]$ as $s$ varies in $S$, hence 
$w_0$ inherits the smoothness of $\lambda_0$, i.e.\ $w_0 \in C^\beta(S)$. 

Now observe that for $\eps > 0$, 
\begin{align*}
 \PP(\|\lambda^*\sigma(g) & -\l_0\|_\infty \le 2 \e)\\ 
& = \PP(\|(\lambda^*-2\|\lambda_0\|_\infty)\sigma(g)+2\|\lambda_0\|_\infty(\sigma(g)-\sigma(w_0))\|_\infty \le 2 \e)\\
& \ge \PP(|\lambda^*-2\|\lambda_0\|_\infty| \le \e)\PP(\|\sigma(g)-\sigma(w_0)\|_\infty \le \e/2\|\lambda_0\|_\infty).
\end{align*}
Since $\l^*$ has a positive, continuous density the first factor on the right is bounded from below by a constant times $\e$. Since the function $\sqrt\sigma$ is Lipschitz by assumption, the second factor is bounded from below by 
$\PP(\|g-w_0\|_\infty \le c\e)$ for a constant $c > 0$. 
By Theorem 3.1 in \cite{vdVvZGamma} we have the lower bound
\[
\PP(\|g-w_0\|_\infty \le \delta_n) \ge e^{-n \delta_n^2},
\]
with $\delta_n$ as specified above. The proof of \eqref{eq: tp1} is now easily completed.

\subsection{Construction of sieves}

Let $\HHH^\ell$ be the RKHS of the GP $g$ with covariance \eqref{eq: g} 
and let $\HHH^\ell_1$ be its unit ball (see \cite{vdVvZRKHS} for background on these notions).
Let $\BBB_1$ be the unit ball in $C[0,1]^d$ relative to the uniform norm. Define
\[
\FF_n = \bigcup_{\lambda\le\lambda_n}\lambda \sigma(\GG_n),
\]
where 
\[
\GG_n = \left[M_n\sqrt{\frac{r_n}{\gamma_n}}\HHH_1^{r_n}+ \eps_n\BBB_1\right]\cup\left[\bigcup_{a \le \gamma_n}(M_n\HHH_1^a) + \eps_n\BBB_1\right],
\]
and $\lambda_n$, $M_n$, $\gamma_n$, $r_n$ and $\eps_n$ are sequences to be determined later. 
In the next two subsections we 
study the metric entropy of the sieves $\FF_n$ and the prior mass of their complements. 

\subsection{Entropy}
Since $\sqrt{\sigma}$ is bounded and Lipschitz we have, for $a,b \in [0,\lambda_n]$, some $c>0$ and $f, g \in \GG_n$, 
\[
\|\sqrt{a\sigma(f)} - \sqrt{b\sigma(g)}\|_\infty \le |\sqrt{a}-\sqrt{b}| + c\sqrt{\lambda_n} \|{f}-{g}\|_\infty.
\]
Since $|\sqrt{a}-\sqrt{b}|\leq \sqrt{|a-b|}$ for $a, b > 0$,  it follows that for $\eps > 0$, 
\[
N(2\eps\sqrt{\lambda_n}, \sqrt{\FF_n}, \|\cdot\|_2) \le N(\eps\sqrt{\l_n}, [0,\l_n], \sqrt{|\cdot|}) N(\eps/c, \GG_n, \|\cdot\|_\infty), 
\]
and hence 
\[
\log N(2\eps\sqrt{\lambda_n}, \sqrt{\FF_n}, \|\cdot\|_2) \lle \log \Big(\frac1\eps\Big) + \log N({\eps}/{c}, {\GG_n}, \|\cdot\|_\infty). 
\]
By formula (5.4) from \cite{vdVvZGamma}, 
\begin{multline*}
\log N(3\eps_n, \GG_n, \|\cdot\|_\infty)
\le K r_n^d \left(\log{\frac{d^{1/4}M_n^{3/2}\sqrt{2\tau r_n}}{\eps_n^{3/2}}}\right)^{1+d}+2\log\frac{2M_n\sqrt{||\mu||}}{\eps_n},
\end{multline*}
for $\|\mu\|$ the total mass of the spectral measure $\mu$, $\tau^2$ the second moment of $\mu$, a constant
 $K>0$, $\gamma_n=\eps_n/(2\tau\sqrt{d} M_n)$, $r_n>A$ for some constant $A>0$, and given that the following 
 relations hold:
\begin{equation}
d^{1/4}M_n^{3/2}\sqrt{2\tau r_n}>2\eps_n^{3/2}, \qquad M_n\sqrt{||\mu||}>\eps_n.
\label{about_entr_1}
\end{equation}
By substituting $\bar\eta_n=\eps_n \sqrt{\lambda_n}$ we get that for some constants $K_1$ and $K_2$,
\[
\log N(2\bar\eta_n, \sqrt{\FF_n}, \|\cdot\|_2) \lle K_1 r_n^d \left(\log \frac{\lambda_n^{3/4} M_n^{3/2} d^{1/4}\sqrt{2\tau r_n}}{\bar\eta_n^{3/2}}\right)^{1+d}+ K_2 \log\frac{\lambda^{1/2}_n M_n}{\bar\eta_n},
\]
when $M_n~>~1$. In terms of $\bar\eta$ the conditions (\ref{about_entr_1}) can be rewritten as
\begin{equation}
\label{about_entr_1_1}
d^{1/4} M_n^{3/2}\lambda_n^{3/4}\sqrt{2\tau r_n}>2\bar\eta_n^{3/2}, \qquad M_n \lambda_n^{1/2} \sqrt{||\mu||}>\bar\eta_n.
\end{equation}
So we conclude that we have the entropy bound
\[
\log N(\bar\eta_n, \sqrt{\FF_n}, \|\cdot\|_2) \lle n \bar\eta_n^2
\]
for sequences $\lambda_n$, $M_n$, $r_n$ and $\bar\eta_n$ satisfying (\ref{about_entr_1_1}) and
\begin{gather}
\begin{split}
\label{about_entr_2}
K_1 r_n^d \left(\log \frac{\lambda_n^{3/4} M_n^{3/2} d^{1/4}\sqrt{2\tau r_n}}{\bar\eta_n^{3/2}}\right)^{1+d}<n\bar\eta_n^2, \quad
K_2 \log\frac{\lambda_n^{1/2} M_n}{\bar\eta_n}<n\bar\eta_n^2.
\end{split}
\end{gather}

\subsection{Remaining mass}
By conditioning we have
\begin{align*}
\PP(\lambda^*\sigma(g) \not \in \FF_n ) & = \int_0^\infty \PP(\l\sigma(g) \not \in \FF_n )p_{\l^*}(\lambda)\,d\lambda\\
& \le\int_0^{\lambda_n} \PP(\l\sigma(g) \not \in \FF_n )p_{\l^*}(\lambda)\,d\lambda + \int_{\lambda_n}^\infty
p_{\l^*}(\lambda)\,d\lambda. 
\end{align*}
By \eqref{eq: flambda} the second term is bounded by a constant times $\exp(-c_0\l_n^\kappa)$. 
For the first term, note that 
for $\lambda \le \lambda_n$ we have
\[
\lambda^{-1}\bigcup_{\lambda' \le \lambda_n}\lambda' \sigma(\GG_n) \supset
\sigma(\GG_n),
\]
hence
$\PP(\lambda\sigma(g) \not \in \FF_n ) \le \PP(g \not \in \GG_n)$.
From (5.3) in \cite{vdVvZGamma} we obtain the bound
\[
\PP(g \not \in \GG_n) \le \frac{K_3 r_n^{p-d+1} e^{-D_2 r_n^d \log^{q}{r_n}}}{\log^{q}{r_n}}+e^{-M_n^2/8},
\]
for some $K_3>0$, $\eps_n<\eps_0$ for a small constant $\eps_0>0$, and $M_n$, $r_n$ and $\eps_n$ satisfying
\begin{equation}
\label{about}
M_n^2>16 K_4 r_n^d (\log(r_n/\eps_n))^{1+d}, \qquad r_n>1,
\end{equation}
where $K_4$ is some large constant.
It follows that $\PP(g \not \in \GG_n)$ is bounded above by a multiple of $\exp{(-L n \widetilde{\eta}_n^2)}$
for a given constant $L$ and $\widetilde{\eta}_n= \lambda_n\eps_n$, provided $M_n$, $r_n$, $\gamma_n$ and $\eps_n$ satisfy
 (\ref{about}) and
\begin{gather}
\begin{split}
\label{about2}
D_2 r_n^d \log^q r_n \ge 2L n \widetilde{\eta}_n^2, \quad
r_n^{p-d+1} \le e^{L n \widetilde{\eta}_n^2}, \quad
M_n^2 \ge 8 L n \widetilde{\eta}_n^2.
\end{split}
\end{gather}
Note that in terms of $\tilde\eta_n$, (\ref{about}) can be rewritten as
\begin{equation}
\label{about1}
M_n^2>16 K_4 r_n^d (\log(r_n\lambda_n/\widetilde{\eta}_n))^{1+d}, \qquad r_n>1.
\end{equation}
We conclude that if (\ref{about1}),(\ref{about2}) holds and 
\begin{equation}
\label{about3}
c_0 \lambda_n^\kappa > L n \widetilde{\eta}_n^2,
\end{equation}
then 
\[
\PP(\l^*\sigma(g \not \in \FF_n )) \lle e^{-L n \widetilde{\eta}_n^2}.
\]

\subsection{Completion of the proof}

In the view of the preceding it only remains to show that $\widetilde{\eta}_n$, $\bar\eta_n$, $r_n$, $M_n>1$ and $\lambda_n$
can be chosen such that relations (\ref{about_entr_1_1}), (\ref{about_entr_2}), (\ref{about2}), (\ref{about1}) and (\ref{about3}) hold.

One can see that it is true for $\widetilde{\eta}_n=\delta_n$ and $\bar\eta_n = \bardelta_n$ described in the theorem, with $r_n$, $M_n$, $\lambda_n$ as follows:
\begin{gather*}
\begin{split}
r_n&=L_2\, n^{\frac{1}{2\beta+d}}(\log n)^\frac{2 k_1}{d}, \\
M_n&=L_3\, n^{\frac{d}{2(2\beta+d)}}(\log n)^{\frac{d+1}{2} + 2k_1}, \\
\lambda_n&=L_4\, n^{\frac{d}{\kappa(2\beta+d)}}(\log n)^\frac{4 k_1}{\kappa}
\end{split}
\end{gather*}
for some large constants $L_2, L_3, L_4 > 0$.

\section{Concluding remarks}
\label{sec: con}

We have shown that the SGCP approach to learning intensity functions proposed 
by \cite{comp_science} enjoys very favorable theoretical properties, provided the 
priors on the hyper parameters are chosen appropriately. 
The result shows there is some flexibility in the construction of the prior. 
The squared exponential GP may be replaced by other smooth stationary processes, 
other link functions may be chosen, 
and there is also a little room in the choice of the priors on the length scale and the multiplicative parameter.
This flexibility is limited, however, and although our result only gives upper bounds on the 
contraction rate, results like those of \cite{Castillo} and 
\cite{jlmr} lead us to believe that one might get sub-optimal performance when deviating 
too much from the conditions that we have imposed.
Strictly speaking the matter is open however and additional research is necessary   to make this belief precise and to describe the exact 
boundaries between good and sub-optimal behaviours.

We expect that a number of generalizations of our results are possible. For instance, 
it should be possible to obtain generalizations to anisotropic smoothness classes 
and priors as 
considered in \cite{anisotropic}, and classes of analytic functions as studied in \cite{vdVvZGamma}.
These generalizations take considerable additional technical work however and 
are therefore not worked out in this paper. We believe they would 
not change the general message of the paper.

\def\cprime{$'$}

\end{document}